\documentclass{article}
\usepackage{epsfig,amsmath,amsthm,amssymb,graphics,verbatim,amsfonts}

\usepackage{makeidx}

\theoremstyle{plain}
\newtheorem{thm}{Theorem}[section]

\newtheorem{prop}[thm]{Proposition}

\theoremstyle{definition}

\def\R{\mathbb{R}}

\def\Z{\mathbb{Z}}
\def\I{\infty}

\newcommand{\be}{\begin{equation}}
\newcommand{\ee}{\end{equation}}
\newcommand{\bea}{\begin{eqnarray}}
\newcommand{\eea}{\end{eqnarray}}
\newcommand{\beann}{\begin{eqnarray*}}
\newcommand{\eeann}{\end{eqnarray*}}
\newcommand{\ra}{\rightarrow}

\begin{document}
 
\date{}
\title{Exploring Parameter Spaces in Dynamical Systems}
\author{Christian Kuehn}

\maketitle

\begin{abstract}
The parameter space of dynamical systems arising in applications is often found to be high-dimensional and difficult to explore. We construct a fast algorithm to numerically analyze ``quantitative features'' of dynamical systems depending on parameters. Using a classical problem from mathematical ecology as an example, we demonstrate how to apply the algorithm to investigate the amplitude of a limit cycle depending on seven parameters. We stress the practical value of the algorithm but we also provide a rigorous error analysis to justify the overall strategy. Our approach turns out to be particularly useful in the case of comparing experimental data to a model defined by differential equations and to investigate whether the equations can approximate the modeled system. 
\end{abstract}

\section{Introduction}

In this paper we investigate the problem of how to explore parameter spaces for dynamical systems. We shall focus on ordinary differential equations (ODE), but remark that our approach works also in different settings such as difference equations, partial differential equations or stochastic dynamics\\

Parametrized dynamical systems in mathematical modelling not only lead to questions about structural features of the model and their dependence on parameters, which is the is central question in applied bifurcation theory, but also to quantitative questions. The main question we are going to address is how a feature, like a limit cycle or a homoclinic orbit, depends on parameters. Some questions often encountered are:
\begin{itemize}
 \item What is the amplitude of a limit cycle?
 \item What is the time required to traverse a part of a homoclinic orbit?
 \item How far does an invariant manifold extend in phase space?
\end{itemize} 
All answers are obviously parameter dependent and are going to have relevant interpretations in the modelling problem under consideration. Drawing quantitative conclusions requires understanding how different parameters influence the feature. Even if the parameter space is of moderate dimension ($\approx 15-25$ parameters) it is often infeasible to numerically explore it systematically by computing a dense grid of parameters values. In many mathematical models no dedicated techniques are used to explore the parameter space. Instead several ``ad-hoc'' methods are used which can potentially yield incorrect conclusions.\\

Most existing techniques for exploring parameter spaces systematically are based on purely statistical methods and in particular on sensitivity analysis. A short exposition of the methods of sensitivity analysis can be found in \cite{SaltelliTarantolaCampolongoRatto}. A detailed explanation of several advanced ideas is given in \cite{SaltelliChanScott}. Our algorithm parallels one of the approaches from sensitivity analysis by using a Monte-Carlo type method. However, the key difference between our approach and known methods is that we take advantage of the vector field ``directly'' in computations to find a more efficient algorithm. In particular, we are not interested in a black-box scheme which works for all mathematical models, but to design a method which uses the structure of the model given by a continuous or discrete-time dynamical system.

\section{Limit Cycles in ODEs - The Spruce Budworm}
\label{sec:budworm}

We use a classical example from mathematical ecology to illustrate our algorithm. The 2-dimensional ODE model proposed in \cite{May} represents the interaction between the spruce budworm $N$ and the leaves of a forest $R$. It is a simplification of a 3-dimensional system for the same ecological situation given in \cite{LudwigJonesHolling}. The paper \cite{LudwigJonesHolling} is partly based on earlier work in \cite{Holling}, whereas re-expositions of the spruce budworm model can be found in many texbooks by now (see e.g. \cite{Yodzis}). Our model reads: 
\bea
\label{eq:model}
\dot{R}&=&p_1 R\left( 1 - \frac{R}{p_3}\right)-p_7 N\nonumber \\
\dot{N}&=&p_2 N\left( 1 - \frac{N}{p_4 R}\right)-p_5\frac{N^2}{p_6^2 R^2+N^2}
\eea
where $R$ is the branch area and $N$ is the budworm density. To keep the notation simple we deviate from the conventions in ecology and simply index the seven parameters $p_1,\ldots,p_7$.\\

\begin{figure}[tbp]
	\centering
		\includegraphics[width=0.7\textwidth]{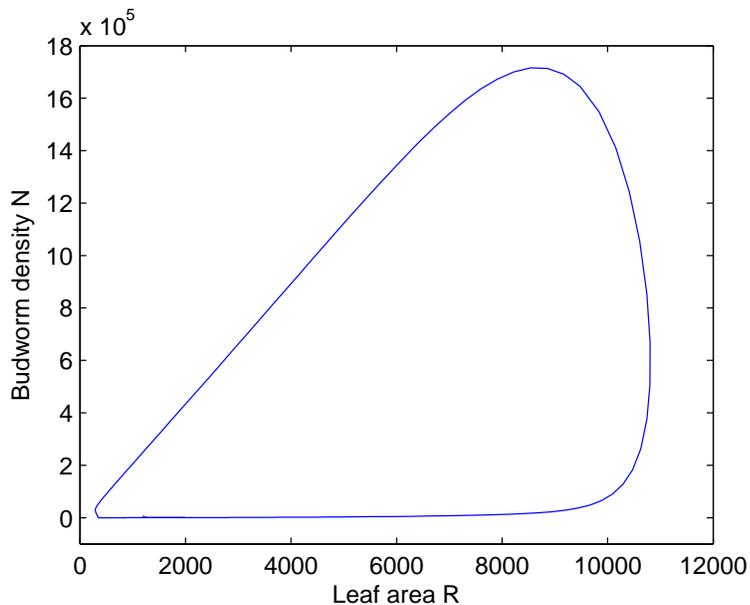}
	\caption{Limit cycle at default parameters given in Table \ref{tab:para}}
	\label{fig:limit_cycle}
\end{figure}

The descriptions, units and default values are summarized in Table \ref{tab:para}. We consider the model only for $R,N\geq0$ representing biologically sensible nonnegative quantities.

\begin{table}[htbp]
	\caption{Parameters for the sample problem}
	\centering
	\label{tab:para}
	\begin{tabular}{|c|l|l|l|}
	\hline
	Parameter & Description & Units & Default Value\\
	\hline\hline
	$p_1$& intrinsic branch growth rate & 1/yr & 0.15 \\
	$p_2$& intrinsic budworm growth rate & 1/yr & 1.6 \\
	$p_3$& maximum branch density & branches/acre & 24000 \\
	$p_4$& maximum budworm density & larvae/branch & 200 \\
	$p_5$& maximum budworm predated & larvae/(acre*year) & 28000 \\
	$p_6$& half of maximum density for predation & larvae/branch & 1.5 \\
  $p_7$& consumption rate of budworm & (branches*acre)/(yr*larvae) & 0.0015 \\
	\hline
	\end{tabular}
\end{table}

We investigate the characteristics of a hyperbolic limit cycle $\gamma\subset\R^2_+$, which represents long-term coexistence for branches (i.e. trees) and the budworm population. A stable limit cycle occurs at the given default parameters (see Figure \ref{fig:limit_cycle}). We observe from equation (\ref{eq:model}) and Table \ref{tab:para} that the variable $R$ is slow compared to the variable $N$. This means that (\ref{eq:model}) is a slow-fast system. Direct calculation of the nullclines shows that the observed limit cycle is a relaxation oscillation.\\


For equation (\ref{eq:model}) we have to explore a 7-dimensional parameter space. We denote the parameter space by $P\subset \R^7_+$ and an element in $P$ by $p=(p_1,\ldots,p_7)$. To indicate the dependence of $\gamma$ on $p$ we shall also refer to the limit cycle as $\gamma(p)$. \\

Biologically relevant questions are to find the maximum values of the coordinates of $\gamma(p)$. These values correspond to the maximum branch area $R_{max}(\gamma(p))$ and maximum budworm density $N_{max}(\gamma(p))$ along a cycle. The maximum budworm density is particularly important as it describes the maximum insect outbreak in the ecosystem. Our main focus is to find the directions in $P$ which correspond to relevant variations of $R_{max}(\gamma(p))$ or $N_{max}(\gamma(p))$. Notice that our aim is to determine parameter dependencies globally as well as locally, i.e. we are interested in the overall behaviour as well as significant variations around points. A further goal is to identify functional dependencies of $R_{max}(\gamma(p))$ or $N_{max}(\gamma(p))$ on subsets of parameters to distinguish between a linear and a nonlinear response.

\section{An Algorithm for Parameter Space Exploration}
\label{sec:algo}

We are going to describe the algorithm in the general setup of ODEs and illustrate all steps with the model given by equation (\ref{eq:model}).

\subsection{The Setup}
\label{sec:setup}
Consider an ordinary differential equation
\be
\label{eq:ode}
\dot{x}=f(x,p) \qquad \text{for $x\in\R^n$, $p\in P\subseteq \R^k$ and $f\in C^1(\R^{n+k})$}
\ee
where $p$ is a vector of parameters $p=(p_1,p_2,\ldots,p_k)$. Suppose that equation (\ref{eq:ode}) has a quantitative property $\gamma(p)$ such as a limit cycle (or e.g. a heteroclinic orbit). By ``quantitative'' we refer to the possibility of measuring an important characteristic of $\gamma(p)$ like cycle amplitude (or e.g. length of a heteroclinic connection). Furthermore assume that we have a compact set $B\subseteq P$ in which the feature is persistent, i.e. we do not have any bifurcation points inside $B$. We shall refer to $B$ as a hyperbolic domain for $\gamma(p)$.\\

Define a function $C:B\rightarrow \R$ given by mapping a parameter vector $p$ to $C(p)$, where $C(p)$ is the quantitative measure of $\gamma(p)$. Assume we have a method to compute $C$ numerically. Notice that $C$ might be very expensive to compute as a single function evaluation usually requires at least one complete numerical solution of equation (\ref{eq:ode}).\\

In our sample problem we take $B=\{p^{default}_i\pm\delta_i\}\subseteq \R^7_+$ where $\delta_i$ are chosen based on the combinations and ranges of parameters we are interested in and to preserve the limit cycle $\gamma(p)$ found in the first quadrant for the default values given in Table \ref{tab:para}. If we compute $\gamma(p)$ numerically we obtain a finite sequence of values $(\gamma(p)_{R,k}, \gamma(p)_{N,k})$, where the indices ``R'' and ``N'' denote the coordinate directions in (\ref{eq:model}). If we are interested in the maximum insect outbreak we define $C(p)$ by
\be
C(p)=\max_{k} \gamma(p)_{N,k}
\ee
Defining $C(p)$ obviously depends on the given problem and the questions of interest in the mathematical model.

\subsection{Discretization}

Assume without loss of generality that $B$ is connected; if this is not the case we can compute for each of the components of $B$ separately. Choose a finite grid $G\subset B$ of spacing $s_i$ in direction $p_i$ defined by intersecting a scaled lattice with $B$ so that 
\be
G=\{ (s_1\cdot \Z,s_2\cdot \Z,\ldots,s_k\cdot \Z)\in \R^k \} \cap B
\ee
Note that we have chosen the lattice to be uniformly spaced in each direction. This is not necessary for any of the following calculations. We should change to a nonuniform grid if there is a priori information about parameter regions available having large variations in $C(p)$. In our model problem we have chosen a uniform spacing in each direction.\\

Now consider the graph formed from $G$ by connecting all points which differ in one coordinate $p_i$ by precisly $s_i$. We also use $G$ to denote this graph. Metrize $G$ by
\be
d(p,q)=\text{minimum number of edges connecting $p$ to $q$}
\ee

\subsection{Relevance Measure}
\label{sec:rel}
The next step is to compare two given elements $p,q\in G$ and estimate when $|C(p)-C(q)|$ is ``large''. To do this we use the given vector field $f$. In the following we shall describe the most basic version of the procedure.\\

\textit{Step 1:} Select a set of points in phase space determined by $\gamma(p)$. This is best illustrated by a concrete example. For our limit cycle $\gamma$ we can choose $k_1$ different parameter values. We compute $\gamma$ for each parameter value and we select $m$ points in phase space on the computed limit cycle. This gives a collection of points in phase space:
\be
S_x:=\{ x_{i} | \quad i=1,2,\ldots,k_2=m k_1 \}
\ee
We choose $k_1$ and $m$ such that the number of elements in $S_x$ is small. In our model problem we have computed $\gamma$ for a set of values chosen uniformly at random from $G$ with $k_1=4$. Then we selected $m=4$ points on each computed limit cycle starting with a random point and uniform spacing along the cycle. This gives a sample representing the different shapes of possible cycles in phase space.\\

\textit{Step 2:} Choose a set points from parameter space
\be
S_p=\{p^1,p^2,\ldots,p^{k_3} \} \subset G
\ee
uniformly at random. Fix a neighbourhood size in $G$, i.e. pick an integer $n_{size}$ such that a neighbourhood of a point $p$ is given by $\Omega_p=\{q\in G|d(q,p)\leq n_{size} \text{ and }q\neq p\}$. For each value $p^j\in S_p$ choose a fixed number of neighbours $k_4$ at random from $\Omega_{p^j}$, call this set $S_{p^j}$; usually we want to pick between $k_4$ between 2 and 4, but any number less than the number of elements in $\Omega_p$ is possible. In preparation for the next step, we mention that we should choose $\#(S_p)=k_3$ such that a summation of $k_4\cdot k_3\cdot k_2$ elements is computationally feasible. \\

\textit{Step 3:} Estimate a measure of relevance. Again we consider the model problem for illustration. Let
\be
\label{eq:rm}
r= \frac{1}{k_2k_3k_4}  \sum_{j=1}^{k_3} \sum_{q\in S_{p^j}} \left| \sum_{i=1}^{k_2} \left(\left\|f(x_{i},p^j)\right\|_2- \left\|f(x_{i},q)\right\|_2 \right)\right|
\ee
So $r$ represents the average variation between two neighbours along the sampled points and hence can be taken as an estimator for the variation of the limit cycle from one point to a neighbouring point. Hence we have found a way to estimate when $|C(p)-C(q)|$ is ``large''. We call $r$ the (average) relevance measure. We can also define 
\be
M(p)=\frac{1}{k_2}\sum_{i=1}^{k_2} \left\|f(x_{i},p)\right\|_2
\ee
which we refer to as the relevance function. We can easily compute the relevance function for different values of $p$ as we only need a small number of evaluations of the vector field if $S_x$ is small. We shall discuss more advanced methods of defining $r$ and $M(p)$ in Section \ref{sec:locexplore}. 

\subsection{Neighbour Comparison}
\label{sec:neigh}

Let $p,q\in G$ be given in parameter space and suppose we have already computed $C(p)$. Now we choose a neighborhood $\Omega_p$ of $p$ in $G$ excluding $p$ itself, e.g. $\Omega_p=\{q\in G|d(p,q)=1\}$. Fix a tolerance $tol\in [0,\I)$, then the rationale for approximating the desired function $C$ is: 
\begin{itemize}
 \item If $q\in \Omega_p$ and $|M(p)-M(q)|\geq tol\cdot r$ $\Rightarrow$ $q$ shows a relevant variation so we compute $C(q)$
 \item If $q\in \Omega_p$ and $|M(p)-M(q)|< tol\cdot r$ $\Rightarrow$ $q$ does not show a relevant variation so set $C(q)=C(p)$
\end{itemize}
Clearly we have two special cases:
\begin{itemize}
 \item $tol=0$: In this case we always compute any point encountered by the algorithm exactly.
 \item $tol=1+\max_{p,q\in G}|M(p)-M(q)|$: In this case we never compute any neighbours of points and just set their values to be the value of the center $p$.  
\end{itemize} 
In our model problem we used the criterion for several values of $tol$ between $0.9$ and $1.1$. The reason for comparing with neighbours will be explained in sections \ref{sec:interpol} and \ref{sec:visual}. 

\subsection{Iteration}
\label{sec:iteration}

We describe next how to run the exploration algorithm to compute $C(p)$. We have two choices to proceed:

\begin{enumerate}
 \item Fix a subset $G_0$ of $G$, for each point $p\in G_0$ compute $C(p)$. Fix a neigbourhood size for each point such that neighbourhoods of points in $G_0$ do not intersect and perform neighbourhood comparisons and computations as described in section \ref{sec:neigh}. We call this method deterministic exploration.
 \item Start with a random initial point $p\in G$. Compute $C(p)$ and mark $p$ as checked. Consider a neighbourhood $\Omega_p$ and perform neighbourhood comparisons and computations as described in section \ref{sec:neigh} for all neighbours which have not been marked yet. Mark all elements in $\Omega_p$. Now choose a new random point out of all unmarked points. We call this method (uniform) random exploration. 
\end{enumerate}

For both methods we can choose a fixed number of points $i_{max}$, which should be explored. Clearly $i_{max}$ should be defined so that $i_{max}\leq \#(G)$. Having reached $i_{max}$ we have obtained the function $C$ on a subset $G^*$ of the grid $G$. Notice that if the grid size shrinks to zero and $i_{max}$ approaches $\#(G)$ the algorithm converges to a full computation of $C$.\\

The most successful algorithms for higher-dimensional numerical problems often use an element of randomization. A classical example is numerical integration. In this area Monte-Carlo techniques are predominant for any problems of dimension $> 4$ (see e.g. \cite{EvansSwartz} for a survey). Hence we expect the deterministic exploration to work only for small number of dimensions and recommend the random exploration.

\subsection{Interpolation}
\label{sec:interpol}  

Given the values of $C$ on $G^*$ we obtain all values in the hyperbolic domain $B$ by interpolation. This constructs $C$ on the full domain $B$. For our model problem we used linear interpolation. For details on basic numerical methods for interpolation consider \cite{Gautschi}. Most algorithms currently in use for multi-dimensional interpolation are based on a Delauney triangulation of given data. Many use the 'quick hull'-algorithm (QHULL, see \cite{BarberDobkinHuhdanpaa},\cite{ORourke}); we used the version of QHULL implemented in MatLab for our sample problem.\\

\begin{figure}[htbp]
	\centering
		\includegraphics[width=1.00\textwidth]{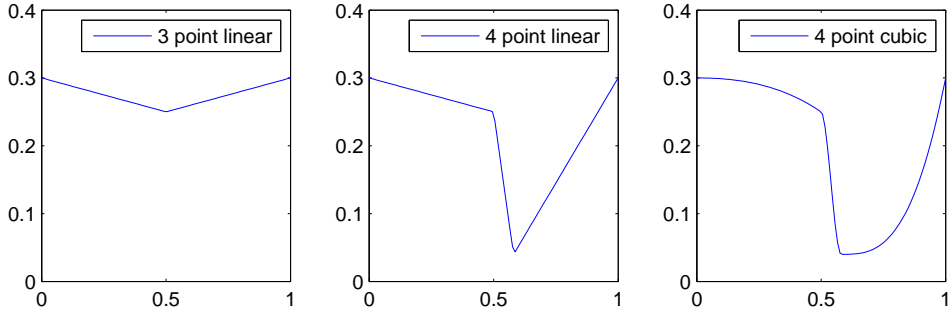}
	\caption{Interpolation for 3 points linear, 4 points linear and 4 points cubic spline}
	\label{fig:interpolation}
\end{figure}

One reason why we have used neighbour comparison as described in section \ref{sec:neigh} is obvious from the one-dimensional toy-example illustrated in Figure \ref{fig:interpolation}. Independent of all previous considerations, just suppose the 1-dimensional interpolation problems on $[0,1]$ are:
\bea
\{(0,0.3),(0.5,0.25), (1,0.3)\} \quad \text{Problem 1}\nonumber\\
\{(0,0.3),(0.5,0.25), (0.58,0.04), (1,0.3)\} \quad \text{Problem 2}\nonumber
\eea

One extra point gives the indication that there is a relevant decrease of the function value to the right of $0.5$ whereas no neighbour comparison would entirely miss this decrease. This means that the local comparion helps us to detect globally important directions of variation of $\gamma(p)$, which could be missed just using point-wise Monte-Carlo exploration without neighbour comparison. It also indicates directions which are of interest near a point. We refer to this feature of our method as local/global approximation.\\

Furthermore we see that other interpolation methods, such as the cubic spline in Figure \ref{fig:interpolation}, work as well, but might introduce a functional dependence due to the interpolation method. Therefore the conservative approach is to use linear interpolation and investigate interesting regions in parameter space more closely once a significant decrease or increase in the function $C$ is found.

\subsection{Results and Visualization}
\label{sec:visual}

To determine the influence of parameters on the structurally stable feature $\gamma(p)$ we have to organize the output of the previously described steps, i.e. we can either work with a subset of values $G^*$ or its interpolation to the hyperbolic domain $B$. Some problems arising are:

\begin{itemize}
 \item Are there ``enough'' points to perform interpolation on $B$ or on some given subset of $B$?
 \item How to visualize the output of high-dimensional data?
\end{itemize}

We shall not discuss the two questions in detail as they form subjects on their own and just focus on a few special issues.\\

The method of neighbour comparison balances the effect of dimension with respect to the number of points necessary for interpolation. Consider our model problem and suppose we fix five parameters on $G^*$ and want to interpolate $C(p)$ over the remaining two parameters, i.e. over a 2-dimensional ``slice'' of $B$. Using just Monte-Carlo simulation without neighbour comparison will inevitably lead to many slices, which do not have any points on them. With neighbour comparison, by possibly increasing the neighbourhood size, we can assure that we have enough points to perform interpolation. Note that higher-dimensional interpolation techniques can also fail if there are too few points. For example, a Delauney triangulation might not be supported by many algorithms for sparse $G^*$ (see \cite{BarberDobkinHuhdanpaa}, \cite{ORourke} and www.qhull.org).\\

For a $k$-dimensional parameter space with $G^*\subset \R^k$ we can fix $k-2$ parameters and plot 3-dimensional output for $C(p)$. This can be achieved by viewing $C(p)$ as a surface over the two variable parameters. We found this method particularly useful coupled to animations, which allow for moving a third parameter and hence produce moving surfaces. During the algorithm we can easily keep track of points, which have been computed and those, which have been set to a local constant by neighbour comparison. Therefore we are able to identify regions and directions containing relevant variations.\\

It should be clear that it is also easy to get estimates of maximum and minimum values over $G^*$ and hence over $B$. For example, in the budworm problem posed in section \ref{sec:budworm} we can now determine for any given region of parameter space an estimate of the maximum insect outbreak.
 
\subsection{Summary of the Algorithm}
The main steps can be summarized as follows:

\begin{enumerate}
 \item Identify the quantitative feature $\gamma(p)$ and its hyperbolic parameter domain $B$.
 \item Define a quantitative measure of $\gamma(p)$ encoded in a function $C:B\rightarrow\R$.
 \item Discretize the space $B$ and define a metric on the discretization $G$.
 \item Sample points from phase space ``representing'' the feature $\gamma(p)$.
 \item Choose a measure (e.g. the norm) to estimate the variation of $\gamma(p)$ from the vector field $f$.
 \item Compute the relevance measure $r$ and setup the relevance function.
 \item Choose discrete or random exploration and a size of neighbourhoods for computed points.
 \item Define a neighbour comparison based on $r$ and the relevance function.
 \item Iteratively compute $C$ on a certain number of points using exploration.
 \item Interpolate (linearly) between the values obtained for $C$.
 \item (Visualize and/or evaluate the interpolated function $C$ to explore parameter dependence.)  
\end{enumerate}

The overview shows that the algorithm works for other dynamical systems like maps, difference equations or partial differential equations.

\section{Advanced Local Exploration}
\label{sec:locexplore}

Instead of using the norm of the vector field, we can also use derivatives in the parameter directions, which leads to the relevance measure:
\be
r_{\partial}:=\frac{1}{k_2k_3k_4k} \sum_{j=1}^{k_3}  \sum_{q\in S_{p^j}} \left|\sum_{i=1}^{k_2} \sum_{l=1}^k\left( \left| \frac{\partial f}{p_l}(x_{i},p^j) \right| -   \left| \frac{\partial f}{p_l}(x_{i},q)\right| \right)\right|
\ee
where we recall that $k$ is the dimension of the parameter space. The corresponding relevance function is given by:
\be
M_{\partial}(p):=\frac{1}{kk_2}\sum_{i=1}^{k_2}   \sum_{l=1}^k \left| \frac{\partial f}{p_l}(x_{i},p) \right| 
\ee

\textit{Remark}: Using the partial derivatives is a technique employed in local sensitivity analysis and has the potential to produce a better measure of relevance. We can use finite-difference approximations to the derivatives if there is no explicit formula of $f(x,p)$ avaiable. $M_{\partial}(p)$ has the drawback that the computational cost increases.

\section{Numerical Results}
\label{sec:numresults} 

Note that since equation (\ref{eq:model}) is a fast-slow system, it is stiff and explicit numerical solvers are going to fail. We used a modified Rosenbrock method of order two with explicitly supplied Jacobian matrix; see \cite{Rosenbrock} for historical reference, \cite{ShampineReichelt} for the implementation and \cite{HairerWannerII} for a modern exposition of the Rosenbrock method. We have simulated the system for the following range of parameters:

\begin{table}[htbp]
	\caption{Ranges of parameters for the sample problem - spruce budworm model}
	\centering
	\label{tab:para2}
	\begin{tabular}{|c|l|l|l|}
	\hline
	Parameter & Default Value & Range & Spacing $s_i$\\
	\hline\hline
	$p_1$& 0.15  & [0.149,0.151] & 0.0001\\
	$p_2$& 1.6   & [1.5,1.7]     & 0.1   \\
	$p_3$& 24000 & [22000,26000] & 250   \\
	$p_4$& 200   & [190,210]     & 2     \\
	$p_5$& 28000 & [24000,32000] & 500   \\
	$p_6$& 1.5   & [1,2]         & 0.1   \\
  $p_7$& 0.0015& [0.001,0.002] & 0.0001\\
	\hline
	\end{tabular}
\end{table}

\begin{figure}
	\centering
		\includegraphics[width=1.00\textwidth]{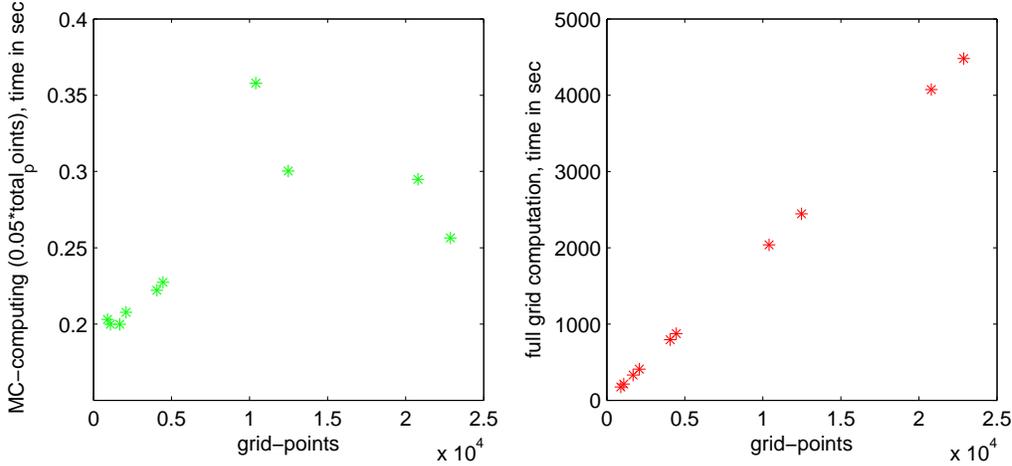}
	\caption{Comparison of computing time for exploration algorithm and full calculation}
	\label{fig:computing_time}
\end{figure}

The main numerical question to pose is, what computational cost do we save by using our algorithm? This is hard to answer precisely as our goal is to provide a method exploring parameter space, which allows applied researchers to detect the main quantitative dependencies in their models. To perform a practical test we have simulated 10 subsets of the full parameter space grid defined in Table \ref{tab:para2}. The results of this simulation are shown in Figure \ref{fig:computing_time}. We have chosen to use random exploration of five percent of the total number of grid points; together with neighbourhood exploration this yields a sufficiently dense set of values to investigate the function $C$.\\

From Figure \ref{fig:computing_time} we observe that the computing time for our algorithm lies between 20 and 40 seconds whereas a full grid computation can take up to 4800 seconds. Now we have to answer the question if we can really extract all the quntitative information we are interested in from our algorithm in comparison to the much slower full grid computation.\\

\begin{figure}
	\centering
		\includegraphics[width=1.00\textwidth]{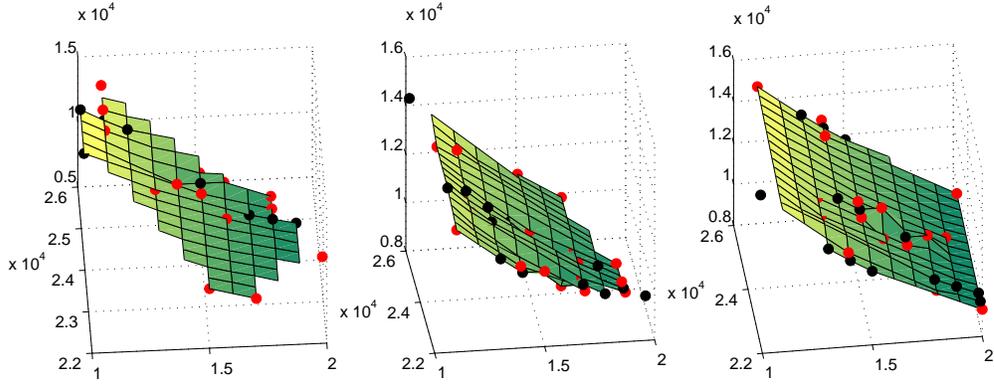}
	\caption{Maximum budworm outbreak as a function of $p_3$ and $p_6$ for $p_5=24000,28000,32000$ - algorithm}
	\label{fig:time_series1}
\end{figure}

Suppose we are interested in the relations between $p_3$, $p_5$ and $p_6$. We simulate the system with $p_3$,$p_5$ and $p_6$ in the ranges given in Table \ref{tab:para2} and use default values for all other parameters as given in Table \ref{tab:para}. We plot the maximum budworm outbreak $N_{max}$ as a function of $p_3$ and $p_6$. Plots for $p_5=24000,28000,32000$ are shown in Figure \ref{fig:time_series1}. The red dots indicate grid points which obtained their values from neighbours. The black dots indicate points which have been fully computed. The computation has been performed with a relevance measure $r$ given in equation (\ref{eq:rm}) and tolerance $tol=1.1$ (see section \ref{sec:neigh}).\\

Just given Figure \ref{fig:time_series1} we can draw several conclusions: Obviously the outbreak is decreasing with increasing value of $p_6$. It is not significantly affected by variations of $p_3$ compared to $p_3$ and $p_6$. We can see that increasing $p_5$ decreases the maximum outbreak by comparing among all three figures. The plots show a weakly nonlinear relationship between $p_3$ and $p_6$. We can also conclude that the quantitative feature varies ``uniformly'' in all three parameter directions, i.e. does not exhibit several nonlinear relationships in our region of interest.\\

\begin{figure}
	\centering
		\includegraphics[width=1.00\textwidth]{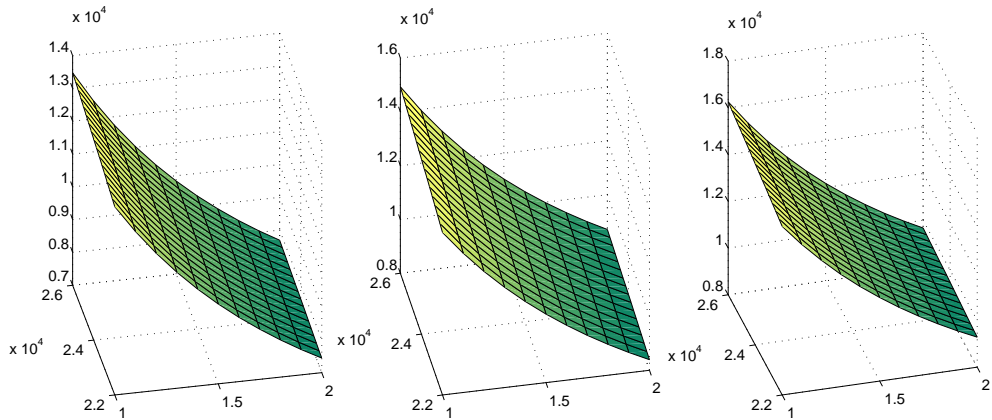}
	\caption{Maximum budworm outbreak as a function of $p_3$ and $p_6$ for $p_5=24000,28000,32000$ - full computation}
	\label{fig:time_series2}
\end{figure}

Now let us compare these findings with the computation on the full grid, which is given in Figure \ref{fig:time_series2}. We see that all our previous conclusions have been verified by the full computation. The only observation we should make is that the first plot in Figure \ref{fig:time_series1} and  Figure \ref{fig:time_series2} do not seem to match well. This results from the fact that our algorithm has produced few values for small $p_6$ and small values of $p_5$ so that this region was not interpolated. We could simply remedy this by extending the current interpolation to the boundary but we decided on the most conservative approach given the computed data.\\

Comparing the computational cost we have computed the two figures on several standard desktop workstations and obtained that it took approximately 10-15 seconds to compute with our algorithm to obtain $N_{max}(p_3,p_5,p_6)$ for the hyperbolic domain $B$ given by the ranges in Table \ref{tab:para2}. This means we have obtained a local/global approximation to $N_{max}$ on $32\cdot16\cdot10=5120$ points. The computation of a single slice of $p_3$ and $p_6$ for a fixed value of $p_5$ took approximately 45-50 seconds on the same machines, yielding a total time for the same amount of points given by $>45 \cdot 16=720$ seconds. This means that in this practical example we have used only a fraction of the full computation
\be
\approx\frac{15}{720}= \frac{1}{48}
\ee 
This 50-fold increase in computation time supplied us with the same conclusions we would have obtained by ``brute force''. Notice that this effect increases with dimension.\\

The computation of a single slice of $p_3$ and $p_6$ using our algorithm is about $1$ second on most current workstations. This means that we can ``move'' with 3-dimensional plots through parameter space in real time and adjust the parameters along directions of interest according to the given output.\\

We have also performed calculations for the full 7-dimensional parameter space and obtained local/global approximations of $N_{max}$ and $R_{max}$. Although this computation takes around 2 hours for most values of $tol\approx 1$, it is still feasible. We remark that the algorithm can also be 'scaled' according to the computational tools available and the desired goals:
\begin{itemize} 
 \item Accurate computations, high computing power, long-waiting time: We decrease $tol$ and select a relatively large fraction of sample points $S_x$.
 \item Coarse global estimates, low computing power, direct feedback: We increase $tol$  and reduce the size of $S_x$.
\end{itemize}
This means that powerful computing environments can be very useful, but are not necessary to extract reliable information from a model.

\section{Error Analysis}
\label{sec:analytic}

The precise error analysis is not the focus of this paper. Nevertheless we shall motivate in this section why we expect an error that behaves like $1/\sqrt{n}$ where $n$ is the number of computed gridpoints. Assume that the hyperbolic domain $B$ is a compact hyperrectangle. First consider the discretizations $G_j\subset B\subset \R^k$ such that the following condition holds:
\bea
\label{eq:cgrid2}
\min_{p\in G_j}\|p-q\|_2\rightarrow 0 \qquad \text{for all $q\in B$ as $j\rightarrow \I$}
\eea
where $\|.\|_2$ denotes the Euclidean norm in $\R^k$. We call the sequence of grids $G_j$ an approximation sequence for $B$ if condition (\ref{eq:cgrid2}) holds. If $f:B\rightarrow \R$ then we define the usual $L^1$ norm as:
\be
\|f\|_{L^1}=\int_B |f(p)|dp
\ee
Since $B$ is a compact hyperbolic domain we conclude that $C:B\rightarrow \R$ is bounded (for the definition of the function $C$, see \ref{sec:setup}). Let $C_A:B\rightarrow \R$ denote the local/global approximation to $C$ obtained by our algorithm for the grid $G_j$. Recall that we have computed $n(j)$ points (depending on $i_{max}$ and the neighboourhood size) on the grid and used linear interpolation to construct $C_A$ on $B$; denote these points by $\{x_k\}_{k=1}^{n(j)}$. Obviously $C_A$ is a random variable depending on the relevance measure $r$ and the grid $G_j$. Hence we should formally write $C_A=C_{A(r,j)}$, but for simplicity of notation we sometimes just write $C_A$. Since $C$ is smooth we can modify $C_A$ on the complement of $G_j$ to be piecewise linear and satisfy $C_A(p)>C(p)$ for all $p\in B$.\\

Let $P$ denote the probability measure on $B$ induced by a uniformly distributed random variable on $B$. Since $B$ is compact we can assume without loss of generality that $\int_Bdx=1$ so that the density of $P$ is identically $1$. Now define
\be
I_A=I_{A(j,r)}=\frac{1}{n(j)}\sum_{k=1}^{n(j)}C_A(x_j)
\ee 
We also set $I=\int_BC(p)dp$. This allows us to consider a standard result from Monte-Carlo integration.

\begin{prop}
\label{prop:emc1}
If $tol=0$ then $C_{ A(r,j) }=C_{A(j)}$, i.e. $C_A$ is indepedent on the random variable $r$. Suppose $n(j)\ra \I$ as $j\ra \I$ then
\be
\frac{I_{A(j)}-I}{\sigma/\sqrt{n(j)}}\stackrel{\sim}{\longrightarrow} Z \qquad \text{ as $j\rightarrow \I$}
\ee
where $Z\sim N(0,1)$, i.e. we have convergence in distribution to a normal distribution with mean 0 and variance 1.
\end{prop}

Remark: Note that $j\rightarrow \I$ implies $n(j)\rightarrow \I$ even if we only increase the number $i_{max}$ linearly. 

\begin{proof}(Sketch, see \cite{EvansSwartz})
If $tol=0$, no neighbour comparison is performed and all elements are computed, so that $C_A$ is independent of $r$ giving the first result.\\

Notice that $Var(I_{A(j)})$ is bounded since $C$ is bounded. Then the central limit theorem implies the desired result.
\end{proof}

For two sequences of random variables $\{X_n\}$, $\{Y_n\}$ we write $X_n=O_P(Y_n)$ if for any $\epsilon>0$ there exist $M_\epsilon>0$, $N_\epsilon>0$ such that
\be
P\left(\left|\frac{X_n}{Y_n}\right|>M_\epsilon\right)<\epsilon \qquad \text{for all $n\geq N_\epsilon$}
\ee 
so that $X_n/Y_n$ is bounded in probability. Therefore we obtain from Proposition \ref{prop:emc1} that:
\be
I_{A(j)}-I=O_{P}\left(\frac{1}{\sqrt{n(j)}}\right)
\ee
Since $O_P(.)\approx O(.)$ and $C_A(p)>C(p)$ for all $p\in B$ we can write
\be
|I_{A(j)}-I|=\|C_A-C\|_{L^1}=O_{P}\left(\frac{1}{\sqrt{n(j)}}\right)\approx O\left(\frac{1}{\sqrt{n(j)}}\right)
\ee
where the last big-oh notation refers to the deterministic asymptotic behaviour on $B$. We can state the previous discussion as a practical result:
\begin{prop}
\label{prop:emc2}
If $tol=0$ then $\|C_A-C\|_{L^1}\approx O\left(\frac{1}{\sqrt{n(j)}}\right)$ for $n(j)$ (respectively $j$) sufficiently large.
\end{prop}

The convergence of the error by an inverse square-root law is obviously very slow. Nevertheless, it is precisely the law we would expect from the use of Monte-Carlo methods. It reflects that we are dealing with a problem without restriction on dimension. Proposition \ref{prop:emc2} justifies our initial claim that we expect an error that behaves like $O(\frac{1}{\sqrt{n}})$. 
  
\section{Conclusions}
\label{sec:last}

We presented a new algorithmic paradigm to explore the quantitative variation of structurally stable features in dynamical systems depending on parameters. We illustrated our algorithm in the case of ODEs and a 2-dimensional model problem and discussed the main steps of discretization, sampling, relevance measures, neighbour comparison, random exploration and interpolation. Numerical results demonstrate the practical value of our method and its efficiency compared to direct computations. Furthermore, error estimates for the random exploration strategy were proven and confirmed the expected Monte-Carlo type behaviour with error of order $O(1/\sqrt{j})$ independent of the dimension of the phase space or parameter space.\\

Despite the detailed error analysis, we stress that our primary goal was not to generate accurate data but to provide practioners with a computationally efficient method. The efficiency derives from a strategy which has been very successful in the analytical description of dynamical systems - we do not compute an explicit solution of an equation but characterize its properties directly from the equations and their structure. We have shown that this general idea can be very useful for the numerical investigation of dynamical systems too.\\

\bibliographystyle{plain}
\bibliography{../my_refs}

\end{document}